\documentclass[reqno]{amsart}

\usepackage{times}
\usepackage{mathabx}
\usepackage{amsmath}
\usepackage{amsthm}
\usepackage{amsfonts}

\newtheorem{thm}{Theorem}[section]

\DeclareMathAlphabet{\mathpzc}{OT1}{pzc}{m}{it}

\numberwithin{equation}{section}

\newcommand{\Wqb}{W_{q,D}}

\newcommand{\R}{\mathbb{R}}

\newcommand{\E}{\mathbb{E}}
\newcommand{\Wq}{\mathbb{W}_q}
\newcommand{\Wqd}{\dot{\mathbb{W}}_q^+}
\newcommand{\Lq}{\mathbb{L}_q}

\newcommand{\Om}{\Omega}

\newcommand{\ve}{\varepsilon}
\newcommand{\rd}{\mathrm{d}}
\newcommand{\divv}{\mathrm{div}_x}

\newcommand{\bqn}{\begin{equation}}
\newcommand{\eqn}{\end{equation}}
\newcommand{\bqnn}{\begin{equation*}}
\newcommand{\eqnn}{\end{equation*}}
\newcommand{\bear}{\begin{eqnarray}} 
\newcommand{\eear}{\end{eqnarray}} 
\newcommand{\bean}{\begin{eqnarray*}} 
\newcommand{\eean}{\end{eqnarray*}} 
\newcommand{\bs}{\begin{split}}
\newcommand{\es}{\end{split}}

\newcommand{\dimens}{\mathrm{dim}}
\newcommand{\codim}{\mathrm{codim}}
\newcommand{\kk}{\mathrm{ker}}
\newcommand{\im}{\mathrm{rg}}

\newcommand{\dhr}{\mathrel{\lhook\joinrel\relbar\kern-.8ex\joinrel\lhook\joinrel\rightarrow}}

\begin{document}

\title[A Nonlocal Nonlinear Reaction-Diffusion Model]{A Note on a Nonlocal Nonlinear Reaction-Diffusion Model}

\author{Christoph Walker}
\email{walker@ifam.uni-hannover.de}
\address{Leibniz Universit\"at Hannover\\ Institut f\" ur Angewandte Mathematik \\ Welfengarten 1 \\ D--30167 Hannover\\ Germany}

\begin{abstract}
We give an application of the Crandall-Rabinowitz theorem on local bifurcation to a system of nonlinear parabolic equations with nonlocal reaction and cross-diffusion terms as well as nonlocal initial conditions. The system arises as steady-state equations of two interacting age-structured populations.
\end{abstract}


\maketitle

\section{Introduction}

In this note we consider coexistence solutions to age-structured population dynamics with diffusion, the main feature of the present work being the inclusion of nonlocal cross-diffusion terms. More precisely, we shall establish positive nontrivial solutions $u=u(a,x),v=v(a,x)$ to the system
\begin{align}
\partial_a u-\divv\big(d_1(\hat{V})\nabla_x u+u\nabla_x d_2(\hat{V})\big) &=-\alpha u^2-\mu_1(\hat{V}) u\ ,\quad a\in (0,a_m)\ ,\quad x\in\Om\ ,  \label{1}\\
\partial_a v-\divv\big(d_3(\hat{U})\nabla_x v+v\nabla_x d_4(\hat{U})\big) &=-\beta v^2-\mu_2(\hat{U}) v\ , \quad a\in (0,a_m)\ ,\quad x\in\Om\ , \label{2} 
\end{align}
for $a\in (0,a_m)$, and $x\in\Om$, subject to the nonlocal initial conditions
\begin{align}
u(0,x)&=\eta U\ , \quad   x\in\Om\ ,\label{3}\\
v(0,x)&=\xi V\ ,  \quad x\in\Om\ ,\label{4}
\end{align}
and Dirichlet boundary conditions
\begin{align}
u(a,x)&=0\ , \quad  a\in (0,a_m)\ ,\quad x\in\partial\Om\ ,\label{5}\\
v(a,x)&=0\ , \quad  a\in(0,a_m)\ ,\quad x\in\partial\Om\ ,\label{6}
\end{align}
where we agree here and in the following upon the notation
\bqn\label{7}
\hat{U}:=\int_0^{a_m}\omega(a) u(a,\cdot)\,\rd a\ ,\qquad U:=\int_0^{a_m} b(a)  u(a,\cdot)\,\rd a
\eqn
for the function $u$ defined on $J:=[0,a_m]$ and analogously for the function $v$. Equations \eqref{1}-\eqref{7} arise naturally as steady-state (i.e. time-independent) equations of two age-structured populations with densities $u$ and $v$, respectively, and maximal age $a_m>0$ living in a (bounded and smooth) domain $\Omega\subset\R^n$, where $a$ is the age and $x$ is the space variable.
The integrals with respect to age in \eqref{1} and \eqref{2} are (weighted) local total populations with a given nonnegative weight function $\omega$. The divergence terms in \eqref{1}, \eqref{2} describe spatial movement with nonlocal coefficients $d_j$. They reflect intrinsic dispersion as well as an increase of dispersive forces by repulsive or attractive interferences with an increase of the other population. We refer to \cite{ShigesadaEtAl} for a derivation of such kind of models (without age-structure). The right hand sides of \eqref{1} and \eqref{2} take into account intra- and inter-specific interactions of the two populations with constants $\alpha, \beta>0$ and functions $\mu_j$ depending nonlocally on the population densities. Creation of new individuals is described by \eqref{3}, \eqref{4} with birth profile $b$ and parameters $\eta$, $\xi$ measuring the intensity of the fertility. We reference to \cite{WebbSpringer} and \cite{WalkerJRAM} for further information on the modeling assumptions.
To avoid unnecessary notational complications, the equations above are stated as a simplified version of more elaborate models, and we remark that the subsequent analysis would not change in any way if one would allow for e.g. different weight functions, different birth rates, or different maximal ages for the two populations. 

In this note we give an extension of previous results \cite{WalkerJRAM, WalkerJFA, WalkerAIHP} being described in more detail in  the next section. Steady-states for a single age-structured populations were investigated e.g. in \cite{DelgadoEtAl2}. We shall also point out that steady-state solutions for two interacting populations when age-structure is neglected, i.e.  variants of the elliptic counterparts of \eqref{1}-\eqref{2}, have attracted considerable interest in the past, see for example \cite{BlatBrown2,DancerTAMS84,DelgadoLGSuarez,Kuto1,KutoYamada,LopezGomezMolinaJDE06,ShiWang}  and the references therein.

\section{Notation and Main Result}

The main features of the equations under consideration are the nonlocalities appearing in the diffusion and reaction terms as well as in the initial conditions.
Similar equations with local reaction terms (i.e. $\mu_1(v)=\alpha_2 v$, $\mu_2(u)=\pm\beta_2 u$ with $\alpha_2, \beta_2>0$) have been investigated in \cite{WalkerJRAM, WalkerJFA} with linear diffusion (i.e. $d_1=d_3= 1$ and $d_2=d_4= 0$) and in
\cite{WalkerAIHP} with a local cross-diffusion term (i.e. $d_1(v)=1+v$, $d_2(v)=v$, $d_3= 1$, and $d_4=0$). In these papers global bifurcation results have been derived with respect to the parameters $\eta$ and $\xi$.
The aim of this note is to show that (local) bifurcation results can be obtained for equations \eqref{1}-\eqref{7} including nonlinear nonlocal diffusion terms. More precisely, we shall provide values for the parameters $\eta$ and $\xi$ for which \eqref{1}-\eqref{7} have {\it coexistence solutions}, i.e. smooth solutions $(u,v)$ with both components nontrivial and positive. Establishing positive steady-state solutions is a first step toward an understanding of (time-dependent) two population dynamics. To this end, we assume throughout the paper that
\bqn\label{8}
d_1\,,\, d_3\,,\, \mu_1\,,\, \mu_2\in C^1(\R)\ ,\quad d_2\,,\, d_4\in C^2(\R)\ ,
\eqn
satisfy
\bqn\label{9}
 d_2(0)=\mu_1(0)=0\ ,
\eqn
respectively,
\bqn\label{10}
d_j(z)\ge\delta\ ,\quad z\ge 0\ ,\quad j=1,3\ ,
\eqn
for some $\delta>0$. For an easier reference in the future we suppose that
\bqn\label{10a}
d_1(0)=1\ .
\eqn
For the weight and the birth functions we assume 
\bqn\label{11}
\omega\,,\, b\in L_\infty^+(J)\ ,\quad b(a)>0\ \text{for $a$ near $a_m$}
\eqn
together with the normalization
\bqn\label{12}
\int_0^{a_m}b(a)e^{-\lambda_1 a}\,\rd a=1\ ,
\eqn
where $\lambda_1>0$ denotes the principal eigenvalue of $-\Delta_x$ on $\Om$ subject to Dirichlet boundary conditions. For technical reasons we introduce $\Lq:=L_q(J,L_q(\Om))$ and the solution space
$$
\Wq:=L_q(J,\Wqb^2)\cap W_q^1(J,L_q)
$$
with $q\in (n+2,\infty)$ fixed, where $L_q:=L_q(\Om)$, and $\Wqb^\kappa:=\Wqb^\kappa(\Om)$ refer to Sobolev-Slobodeckii spaces including Dirichlet boundary conditions if meaningful, i.e. if $\kappa>1/q$. 
We let $\Wq^+$ denote the positive cone of $\Wq$ and set $\Wqd:=\Wq^+\setminus\{0\}$. Recall the embedding 
\bqn\label{emb}
\Wq\hookrightarrow C^{1-1/q-\vartheta}\big([0,a_m],\Wqb^{2\vartheta}\big)\ ,\quad 0\le \vartheta\le 1-1/q\ ,
\eqn
which holds for $\vartheta=1-1/q$ due to \cite[III.Thm.4.10.2]{LQPP} and otherwise by the interpolation inequality \cite[I.Thm.2.11.1]{LQPP}. In particular, the trace $\gamma_0u:=u(0)$ defines a bounded linear operator $\gamma_0\in\mathcal{L}(\Wq,\Wqb^{2-2/q})$. Also recall (e.g. from \cite{LQPP}) that $A\in L_\infty(J,\mathcal{L}(\Wqb^2,L_q))$ is said to have {\it maximal $L_q$-regularity} provided that the operator $$(\partial_a+A,\gamma_0)\in\mathcal{L}(\Wq,\Lq\times\Wqb^{2-2/q})$$ has a bounded inverse, where $(A\phi)(a):=A(a)\phi(a)$ for $a\in J$ and $\phi\in \Wq$. We also note the compact embeddings
\bqn\label{emb2}
\Wqb^2\dhr \Wqb^{2-2/q}\dhr C^1(\bar{\Om})
\eqn
as $q>n+2$, and that the interior of the positive cone $\Wqb^{2-2/q, +}$ of $\Wqb^{2-2/q}$, denoted by $\mathrm{int}(\Wqb^{2-2/q, +})$, is nonempty. Finally, recall that an operator $Z\in\mathcal{L}(\Wqb^{2-2/q}):=\mathcal{L}(\Wqb^{2-2/q},\Wqb^{2-2/q})$ is {\it strongly positive} if $Z\phi\in \mathrm{int}(\Wqb^{2-2/q, +})$ for~$\phi\in \Wqb^{2-2/q, +}\setminus\{0\}$.
\\

On taking $v\equiv 0$, our assumptions imply that we get from \eqref{1}-\eqref{7} a reduced problem with a nonlocal initial condition of the form
\bqn\label{Aee1}
\partial_a u-\Delta_x u=-\alpha u^2\ ,\quad u(0)=\eta U\ ,
\eqn
subject to Dirichlet boundary conditions, which has been studied in a previous work \cite[Thm.2.1]{WalkerJRAM}. In fact, in view of the normalization \eqref{12}, this problem admits no solution in $\Wqd$ if $\eta\le 1$, while for each $\eta>1$ there is a unique solution $u_\eta\in\Wqd$ depending smoothly on $\eta$ and $\|u_\eta\|_{\Wq}\rightarrow\infty$ as $\eta\rightarrow\infty$. Parabolic regularity theory implies that $u_\eta$ is smooth with respect to both $a\in J$ and $x\in\Om$. 

In the following we shall consider the situation that $\eta>1$ is given. We then regard $\xi$ as a bifurcation parameter and write $(\xi,u,v)$ for solutions to \eqref{1}-\eqref{7}. The considerations above guarantee the existence of a semi-trivial branch of solutions
$$
\mathfrak{T}_0:=\{(\xi,u_\eta,0)\,;\, \xi\ge 0\}\ .
$$
On applying the famous Crandall-Rabinowitz theorem on local bifurcation, we derive the existence of a local branch of coexistence solutions bifurcating from the semi-trivial branch $\mathfrak{T}_0$. More precisely, we have:

\begin{thm}\label{T}
Given $\eta>1$, there exists $\xi_0:=\xi_0(\eta)>0$ such that a local curve $\mathfrak{T}$ in $\R^+\times\Wqd\times\Wqd$ of (smooth) coexistence solutions $(\xi,u,v)$ to \eqref{1}-\eqref{7} emanates from $(\xi_0,u_{\eta},0)\in\mathfrak{T}_0$. This bifurcation point is unique, i.e. there is no other bifurcation point on $\mathfrak{T}_0$ to coexistence solutions.
\end{thm}

As equations \eqref{1}-\eqref{7} are symmetric in $u$ and $v$, one may interchange the role of $\xi$ and $\eta$, of course. We remark that establishing a bifurcation from the other semi-trivial branch $\{(\xi,0,v_\xi);\xi>1\}$ when regarding $\xi$ as bifurcation parameter (if $\mu_2(0)=0$ and $d_4(0)=0$) does not seem to be obvious for the present situation. However, for the case of linear diffusion, i.e. if $d_1=d_3=1$ and $d_2=d_4=0$, one can show global bifurcation results for \eqref{1}-\eqref{7} along the lines of \cite{WalkerJRAM, WalkerJFA} for nonlocal reaction terms as well.

\section{Proof of Theorem~\ref{T}}

The remainder is dedicated to the proof of Theorem~\ref{T}, which is a consequence of the bifurcation result of Crandall $\&$ Rabinowitz \cite[Thm.1.7]{CrandallRabinowitz}. Let $\eta>1$ be fixed. We write solutions to \eqref{1}-\eqref{7} in the form $(\xi,u,v)=(\xi,u_\eta-w,v)$, which we then obtain as the zeros $(\xi,w,v)$ of the function $$F:\R\times\Wq\times\Wq\longrightarrow \Lq\times\Lq\times\Wqb^{2-2/q}\times\Wqb^{2-2/q}\ ,$$ 
where 
$$
F(\xi,w,v):=\left(\begin{array}{c} 
\partial_a w-\divv\big(d_1(\hat{V})\nabla_x w+ w\nabla_x d_2(\hat{V})\big)-\divv\big((1-d_1(\hat{V}))\nabla_x u_\eta 
- u_\eta\nabla_x d_2(\hat{V})\big)\\
\qquad\qquad\qquad\qquad\qquad\qquad\qquad-\mu_1(\hat{V})(u_\eta-w)+2\alpha u_\eta w-\alpha w^2\\
\partial_a v-\divv\big(d_3\big(\hat{U}_\eta-\hat{W}\big)\nabla_x v +v \nabla_x d_4\big(\hat{U}_\eta-\hat{W}\big)\big) +\beta v^2+\mu_2\big(\hat{U}_\eta-\hat{W}\big)v\\
w(0)-\eta W\\
v(0)-\xi V
 \end{array}\right)
$$
for $(\xi,w,v)\in \R\times\Wq\times\Wq$.
Owing to \eqref{9} and \eqref{10a}, the Frech\'{e}t derivatives at $(\xi,w,v)=(\xi,0,0)$ are given by
$$
F_{(w,v)}(\xi,0,0)[\phi,\psi]=\left(\begin{array}{c} 
\partial_a \phi-\Delta_x\phi +\divv\big(d_1'(0)\hat{\Psi}\nabla_x u_\eta +u_\eta d_2'(0)\nabla_x\hat{\Psi}\big) -\mu_1'(0)\hat{\Psi}u_\eta +2\alpha u_\eta \phi\\
\partial_a \psi-\divv\big(d_3\big(\hat{U}_\eta\big)\nabla_x \psi +\psi \nabla_x d_4\big(\hat{U}_\eta\big)\big)+\mu_2\big(\hat{U}_\eta\big)\psi\\
\phi(0)-\eta \Phi\\
\psi(0)-\xi \Psi
 \end{array}\right)
$$
with dashes denoting derivatives, and
$$
F_{\xi,(w,v)}(\xi,0,0)[\phi,\psi]=\left(\begin{array}{c} 
0\\
0\\
0\\
- \Psi
 \end{array}\right)\ ,
$$
where we use here and in the following  notation \eqref{7} for $(\phi,\psi)\in\Wq\times\Wq$.\\

{\bf (i)} We shall show that $F_{(w,v)}(\xi,0,0)$ is a Fredholm operator of index zero for $\xi>0$. To this end we introduce operators $A_1\in C^{1/2}(J,\mathcal{L}(\Wqb^2,L_q))$ and $A_3\in\mathcal{L}(\Wqb^2,L_q)$ as
\begin{align*}
A_1(a)\phi &:=-\Delta_x\phi+2\alpha u_\eta(a) \phi\ ,\\
A_3\psi &:=-\divv\big(d_3\big(\hat{U}_\eta\big)\nabla_x \psi +\psi \nabla_x d_4\big(\hat{U}_\eta\big)\big)+\mu_2\big(\hat{U}_\eta\big)\psi\  ,
\end{align*}
for $a\in J$ and $\phi, \psi\in \Wqb^2$, the regularity being implied by \eqref{emb} as $u_\eta\in\Wq$. It follows from \eqref{10} and \cite[III.Sect.4]{LQPP} (in particular, see III.Ex.4.7.3d), III.Thm.4.8.7, and III.Thm.4.10.10 of \cite{LQPP}) that both $A_1$ and $A_3$ have maximal $L_q$-regularity.
We also introduce $A_2\in \mathcal{L}(\Wq,\Lq)$ by
$$
(A_2\psi)(a) :=\divv\big(d_1'(0)\hat{\Psi}\nabla_x u_\eta (a)+u_\eta(a) d_2'(0)\nabla_x\hat{\Psi}\big) -\mu_1'(0)\hat{\Psi}u_\eta(a)\ , \qquad a\in J\ ,\quad \psi\in \Wq\ ,
$$
and set, for $(\phi,\psi)\in\Wq\times\Wq$,
$$
\mathbb{A}(\phi,\psi):=\left(\begin{array}{c}A_1\phi +A_2\psi\\ A_3\psi\end{array}\right)
$$
using the convention $(A_1\phi)(a)=A_1(a)\phi(a)$ and $(A_3\psi)(a)=A_3\psi(a)$ for $a\in J$.
On letting 
$$
\gamma_0(\phi,\psi):=\left(\begin{array}{c}\phi(0)\\ \psi(0)\end{array}\right)\ ,\quad (\phi,\psi)\in\E_1:=\Wq\times\Wq\ ,
$$  
and 
$$
\E_0:=\Lq\times\Lq\ ,\quad E_\varsigma:=\Wqb^{2-2/q}\times\Wqb^{2-2/q}\ ,
$$
it is readily seen from the triangular structure of $\mathbb{A}$ that the operator 
$
\big(\partial_a+\mathbb{A},\gamma_0\big)\in \mathcal{L}(\E_1,\E_0\times E_\varsigma)
$
has a bounded inverse, say, $$T:=\big(\partial_a+\mathbb{A},\gamma_0\big)^{-1}\in\mathcal{L}(\E_0\times E_\varsigma,\E_1)\ .$$ 
Introducing $\ell[\xi]\in\mathcal{L}(\E_1,E_1)$ for $E_1:=\Wqb^2\times\Wqb^2$ by
$$
\ell[\xi](\phi,\psi):=\left(\begin{array}{c}\eta\Phi\\ \xi\Psi\end{array}\right)\ ,\quad (\phi,\psi)\in\E_1\ ,
$$
we obtain from a straightforward modification of \cite[Lem.2.1]{WalkerSIMA} (with $E_0:=L_q\times L_q$ therein) that 
$$
L[\xi]:= F_{(w,v)}(\xi,0,0)=\big(\partial_a+\mathbb{A}, \gamma_0 -\ell[\xi]\big)\in  \mathcal{L}\big(\E_1,\E_0\times E_\varsigma\big)
$$ 
is indeed a Fredholm operator of index zero. In fact, defining
$$
Q[\xi]w:=\ell[\xi]( T(0,w))\ ,\quad w\in E_\varsigma\ ,
$$
and observing that $Q[\xi]\in \mathcal{L}(E_\varsigma,E_1)$ is a compact operator on $E_\varsigma$ into itself due to the compact embedding \mbox{$E_1\hookrightarrow E_\varsigma$} implied by \eqref{emb2}, we have
\begin{align*}
\kk(L[\xi])&=\big\{T(0,w)\,;\ w\in\kk(1-Q[\xi])\big\}\ ,\\
\im(L[\xi])&=\big\{(f,h)\in \E_0\times E_\varsigma\,;\, h+\ell[\xi](T(f, 0))\in \im(1-Q[\xi])\big\}\ ,
\end{align*}
both spaces being closed with
\bqn\label{CR1}
\dimens(\kk(L[\xi]))=\codim(\im(L[\xi]))=\dimens(\kk(1-Q[\xi]))<\infty\ .
\eqn

{\bf (ii)} We now choose $\xi_0>0$ such that $\kk(L[\xi_0])$ is one-dimensional. First observe that, owing to the parabolic maximum principle (e.g. see \cite[Cor.13.6]{DanersKochMedina}), the semigroup $\{e^{-aA_3} ; a\ge 0\}$ on $L_q$ generated by $-A_3$ is such that $e^{-aA_3}\in\mathcal{L}(\Wqb^{2-2/q})$ is strongly positive for $a>0$. Hence, \eqref{11} and standard regularity effects of analytic semigroups \cite{LQPP} ensure that
$$
H:=\int_0^{a_m} b(a)\, e^{-aA_3}\,\rd a\in \mathcal{L}(\Wqb^{2-2/q},\Wqb^2)
$$
defines a strongly positive and compact operator on $\Wqb^{2-2/q}$ into itself by \eqref{emb2}. Its spectral radius $r(H)>0$ is thus a simple eigenvalue of $H$ with a corresponding eigenfunction $\Psi_0\in \mathrm{int}(\Wqb^{2-2/q,+})$ due to the Krein-Rutman theorem \cite[Thm.3.2]{AmannSIAMReview}. We then put
$$
\xi_0:=\xi_0(\eta):=\frac{1}{r(H)}
$$
and obtain $\kk(1-\xi_0 H)=\mathrm{span}\{\Psi_0\}$. Let us also observe that $A_1\in C^{1/2}(J,\mathcal{L}(\Wqb^2,L_q))$ generates a parabolic evolution operator $U_1(a,\sigma)$ on $L_q$ by \cite[Cor.4.4.2]{LQPP}. The same arguments as above ensure that
$$
G_1:=\int_0^{a_m} b(a)\, U_1(a,0)\,\rd a\in \mathcal{L}(\Wqb^{2-2/q},\Wqb^2)\ ,
$$
that $G_1$ is compact on $\Wqb^{2-2/q}$ into itself by \eqref{emb2}, and that $G_1$ is strongly positive by \eqref{11} since this is true for each $U_1(a,0)$, $a\in (0,a_m]$, see \cite[Cor.13.6]{DanersKochMedina}. Next we claim that $1-\eta G_1$ is invertible. Indeed, owing to the positivity of $u_\eta$, the evolution operator $U_2(a,\sigma)$ on $L_q$ generated by $A_1-\alpha u_\eta=-\Delta_x +\alpha u_\eta$ (subject to Dirichlet boundary conditions) is such that $U_2(a,0)- U_1(a,0)$ is strongly positive on $\Wqb^{2-2/q}$ for $a\in (0,a_m]$ by the parabolic maximum principle, whence $G_2-G_1$ is strongly positive due to \eqref{11}, where
$$
G_2:=\int_0^{a_m} b(a)\, U_2(a,0)\,\rd a \ .
$$
Therefore, $r(G_2)>r(G_1)$ by the Krein-Rutman theorem \cite[Thm.3.2]{AmannSIAMReview}. Recall then that \eqref{Aee1} implies $u_\eta(a)=U_2(a,0)u_\eta(0)$ with $u_\eta(0)=\eta G_2 u_\eta(0)$ and so $r(\eta G_2)=1$ again by the Krein-Rutman theorem since $u_\eta(0)$ is strictly positive. Consequently, $1-\eta G_1$ is invertible.
Consider now the equation  $L[\xi_0](\phi,\psi)=0$ for $(\phi,\psi)\in\E_1$, that is,
\bqn\label{zz}
\partial_a\psi+A_3\psi=0 \ ,\quad a\in (0,a_m]\ ,\qquad\psi(0)=\xi_0\Psi\ ,
\eqn
and 
\bqn\label{zzz}
\partial_a\phi+A_1(a)\phi=-(A_2\psi)(a) \ ,\quad a\in (0,a_m]\ ,\qquad \phi(0)=\eta\Phi\ .
\eqn
From \eqref{zz} we deduce $\psi(a)=e^{-aA_3}\psi(0)$ with $\psi(0)=\xi_0 H \psi(0)$, whence $\psi(0)=\zeta\Psi_0$ for some $\zeta\in\R$ since $1$ is a simple eigenvalue of $\xi_0H$ with eigenvector $\Psi_0$, and so $\psi=\zeta\psi_*$ with $\psi_*(a):=e^{-aA_3}\Psi_0$. From the first part of \eqref{zzz} it then follows that 
$$
\phi(a)=U_1(a,0)\phi(0)-\zeta\int_0^a U_1(a,\sigma)\, (A_2 \psi_*)(\sigma)\,\rd \sigma\ ,\quad a\in J\ ,
$$
which, when plugged into the initial condition $\phi(0)=\eta\Phi$, yields $\phi(0)=\zeta \Phi_0$, where
$$
\Phi_0:=-\eta \big(1-\eta G_1\big)^{-1}\int_0^{a_m} b(a)\int_0^a U_1(a,\sigma)\, (A_2\psi_*)(\sigma) \,\rd \sigma\,\rd a\ .
$$
Thus, setting
$$
\phi_*(a):=U_1(a,0)\Phi_0-\int_0^a U_1(a,\sigma)\, (A_2 \psi_* )(\sigma) \,\rd \sigma\ ,\quad a\in J\ ,
$$
we conclude $(\phi,\psi)=\zeta (\phi_*,\psi_*)$ and then
\bqn\label{CR2}
\kk(L[\xi_0])=\mathrm{span}\{(\phi_*,\psi_*)\}\ .
\eqn

{\bf (iii)} Next, we check the transversality condition of \cite{CrandallRabinowitz}, that is, we show that
\bqn\label{CR3}
F_{\xi,(w,v)}(\xi_0,0,0)[\phi_*,\psi_*]\not\in \mathrm{rg}\big(L[\xi_0]\big) \ .
\eqn
Supposing otherwise, there is $\psi\in\Wq$ such that
$$
\partial_a\psi + A_3\psi=0\ ,\quad a\in (0,a_m]\ ,\qquad \psi(0)-\xi_0\Psi=-\Psi_*\ ,
$$
and we easily derive $(1-\xi_0 H)\psi(0)=-\Psi_*$. Since $$\Psi_0\in \mathrm{int}(\Wqb^{2-2/q,+})\cap \kk(1-\xi_0 H)\ ,$$ we may choose $\kappa>0$ sufficiently large such that $$P:=\kappa\Psi_0-\psi(0)\in \mathrm{int}(\Wqb^{2-2/q,+})\ , \quad (1-\xi_0 H)P=\Psi_*\ .$$ 
Because $r(\xi_0 H)=1$ and $\Psi_*\in\Wqb^{2-2/q,+}\setminus\{0\}$, this contradicts the fact that this last equations has no positive solution $P$ according to \cite[Thm.3.2]{AmannSIAMReview}, and we conclude \eqref{CR3}. \\

{\bf (iv)} Gathering \eqref{CR1}, \eqref{CR2}, and \eqref{CR3} we are now in a position to apply \cite[Thm1.7]{CrandallRabinowitz} and deduce that the nontrivial zeros of $F$ close to the bifurcation point $(\xi_0,0,0)$ lie on a continuous curve. More precisely and applied to \eqref{1}-\eqref{7}, we obtain some $\ve_0>0$ and continuous functions $\xi:(-\ve_0,\ve_0)\rightarrow \R$ and $\zeta_j:(-\ve_0,\ve_0)\rightarrow\Wq$ with $\xi(0)=\xi_0$ and $\zeta_j(0)=0$, such that the nontrivial solutions $(\xi,u,v)$ to \eqref{1}-\eqref{7} lie on the curve
$$
\{(\xi(\ve),u_\eta-\ve\phi_*-\ve\zeta_1(\ve), \ve\psi_*+\ve\zeta_2(\ve))\,;\, -\ve_0<\ve<\ve_0\}\ \subset\R\times\Wq\times\Wq\ ,
$$
which bifurcates from $(\xi_0,u_\eta,0)$. Since $u_\eta (0)$ and $\psi_*(0)=\Psi_0$ belong to $\mathrm{int}(\Wqb^{2-2/q,+})$, it follows from \eqref{emb} that, for $\ve>0$ sufficiently small, $u^\ve:=u_\eta-\ve\phi_*-\ve\zeta_1(\ve)$ and $v^\ve:=\ve\psi_*+\ve\zeta_2(\ve)$ both have nontrivial initial values 
$$
u^\ve(0)=u_\eta(0)-\ve\Phi_0-\ve\gamma_0\zeta_1(\ve)\in\Wqb^{2-2/q,+}
$$ 
and 
$$ 
v^\ve(0)=\ve\Psi_0+\ve\gamma_0\zeta_2(\ve)\in\Wqb^{2-2/q,+}\ ,
$$ 
respectively. Therefore, making $\ve_0$ smaller, if necessary, we derive from \eqref{1}, \eqref{2}, and the parabolic maximum principle that
$$
\mathfrak{T}:=\{(\xi(\ve),u_\eta-\ve\phi_*-\ve\zeta_1(\ve), \ve\psi_*+\ve\zeta_2(\ve))\,,\, 0<\ve<\ve_0\}
$$
is a curve of solutions to \eqref{1}-\eqref{7} in $\R^+\times\Wqd\times\Wqd$ bifurcating from $(\xi_0,u_\eta,0)\in\mathfrak{T}_0$.\\

{\bf (v)} Finally, if $(\xi,u_\eta,0)\in\mathfrak{T}_0$ is any bifurcation point to positive coexistence solutions of \eqref{1}-\eqref{7}, there is a sequence $(\xi_j,u_j,v_j)$ in $\R^+\times\Wqd\times\Wqd$ of solutions to \eqref{1}-\eqref{7} converging to this point. Since solutions to \eqref{1}-\eqref{7} are smooth by standard parabolic regularity theory, it is not difficult to see that a subsequence of $\psi_j:=v_j/\|v_j\|_{\Wq}$ converges to some $\psi\in\Wqd$ satisfying 
$$
\partial_a\psi+A_3\psi=0 \ ,\quad a\in (0,a_m]\ ,\qquad\psi(0)=\xi\Psi\ ,
$$
what readily implies $\xi=\xi_0$ since $\psi(0)=\xi H \psi(0)$, i.e. $r(\xi H)=1$ by the Krein-Rutman theorem. Thus ($\xi_0,u_\eta,0)$ is the only bifurcation point on $\mathfrak{T}_0$. This proves Theorem~\ref{T}.


\begin{thebibliography}{99}

\bibitem{AmannSIAMReview}
H. Amann. \textit{Fixed point equations and nonlinear eigenvalue problems in ordered Banach spaces.}  SIAM Rev. {\bf 18} (1976), no. 4, 620-709.

\bibitem{LQPP}
H. Amann. \textit{Linear and Quasilinear Parabolic Problems,
{Volume} {I}: Abstract Linear Theory.} Birkh\"auser 1995.

\bibitem{BlatBrown2}
J. Blat, K.J. Brown. \textit{Global bifurcation of positive solutions in some systems of elliptic equations.} SIAM J. Math. Anal. {\bf 17} (1986), 1339-1353.

\bibitem{CrandallRabinowitz}
M.C. Crandall, P.H. Rabinowitz. {\it Bifurcation from simple eigenvalues.} J. Functional Analysis {\bf 8} (1971), 321-340.

\bibitem{DancerTAMS84}
E.N. Dancer. {\em On positive solutions of some pairs of differential equations.} Trans. Amer. Math. Soc. {\bf 284} (1984), 729-743.

\bibitem{DanersKochMedina}
D. Daners, P. Koch-Medina. {\it Abstract Evolution Equations, Periodic Problems, and Applications.}
Pitman Res. Notes Math. Ser., {\bf 279}, Longman, Harlow 1992. 

\bibitem{DelgadoLGSuarez}
M. Delgado, J.L\'opez-G\'omez, A. Su\'arez. {\it On the symbiotic Lotka-Volterra model with diffusion and transport effects.} J. Differential Equations {\bf 160} (2000) 175-262.

\bibitem{DelgadoEtAl2}
M. Delgado, M. Molina-Becerra, A. Su\'arez. {\it Nonlinear age-dependent diffusive equations: A bifurcation approach.}
J. Differential Equations {\bf 244} (2008), 2133-2155.

\bibitem{Kuto1}
K. Kuto. {\it Stability of steady-state solutions to a prey-predator system with cross-diffusion.} J. Differential Equations {\bf 197} (2004) 293-314. 

\bibitem{KutoYamada}
K. Kuto, Y. Yamada. {\it Multiple coexistence states for a prey-predator system with cross-diffusion.} J. Differential Equations {\bf 197} (2004) 315-348. 

\bibitem{LopezGomezMolinaJDE06}
J. L\'opez-G\'omez, M. Molina-Meyer. {\em Superlinear indefinite systems: beyond Lotka-Volterra models.} J. Differential Equations {\bf 221} (2006), 343-411.

\bibitem{ShigesadaEtAl}
N. Shigesada, K. Kawasaki, E. Teramoto. {\it Spatial segregation of interacting species.} J.
Theoretical Biology, {\bf 79} (1979), 83-89.

\bibitem{ShiWang}
J. Shi, X. Wang. {\it On global bifurcation for quasilinear elliptic systems on bounded domains.}
J. Differential Equations {\bf 246} (2009), no. 7, 2788-2812. 

\bibitem{WalkerSIMA}
Ch. Walker. {\em Positive equilibrium solutions for age and spatially structured population models.}
SIAM J. Math. Anal. {\bf 41} (2009), 1366-1387.

\bibitem{WalkerJRAM}
Ch. Walker. {\em On positive solutions of some system of reaction-diffusion equations with nonlocal~initial conditions.}
To appear in: J. Reine Angew. Math. (arXiv:1003.4698[math.AP]).


\bibitem{WalkerJFA}
Ch. Walker. {\em On nonlocal parabolic steady-state equations of cooperative or competing systems.} Preprint (2010) (arXiv:1008.3125[math.AP]).

\bibitem{WalkerAIHP}
Ch. Walker. {\em Positive solutions of some parabolic system with cross-diffusion and nonlocal initial conditions.} Preprint (2010) (arXiv:1011.3334[math.AP]).


\bibitem{WebbSpringer}
G.F. Webb. {\em Population models structured by age, size, and spatial position.} In: P. Magal, S. Ruan (eds.) {\em Structured Population Models in Biology and Epidemiology.} Lecture Notes in Mathematics, Vol. 1936. Springer, Berlin, 2008.

\end{thebibliography}
\end{document}